\documentclass[12pt]{amsart}
\usepackage{amsmath,amsfonts,amssymb}

\def\le{\leqslant}
\def\ge{\geqslant}

\def\vt{\vartheta}

\renewcommand{\mod}{\mathop{\rm{mod}}}
\def\qb{\mathbb{Q}}

\numberwithin{equation}{section}

\newtheorem{cor}{COROLLARY}[section]
\newtheorem{prop}[cor]{PROPOSITION}

\theoremstyle{definition}

\newtheorem{theorem}{THEOREM}

\begin{document}
\title{\bf A Note on Dirichlet $L$-functions}
\author[Friedlander]{J.B. Friedlander$^*$}
\thanks{$^*$\ Supported in part by NSERC grant A5123}
\author[Iwaniec]{H. Iwaniec$^{**}$}
\thanks{$^{**}$\ Supported in part by NSF grant DMS-1406981}

\maketitle

{\bf Abstract:} 
We study the relation between the size of $L(1,\chi)$ and the width of the 
zero-free interval to the left of that point.

\section{\bf Introduction}

Let $\chi (\mod D)$ be a real, primitive character of conductor $D$ and 
$$
L(s,\chi)=\sum_1^{\infty}\chi (n)n^{-s} 
$$ 
its Dirichlet $L$-function. We are interested in the size of $L(1,\chi)$. The 
Generalized Riemann Hypothsis implies 
\begin{equation}\label{eq:1.1}
(\log\log D)^{-1}\ll L(1,\chi) \ll  \log\log D\ .
\end{equation} 
Unconditionally, we have the easy upper bound 
$$
L(1,\chi) \ll  \log D\ ,
$$
but what one might expect to be the corresponding lower bound, 
\begin{equation}\label{eq:1.2}
L(1,\chi) \gg  (\log D)^{-1}\ ,
\end{equation} 
is not yet known to hold in general. We know that (1.2) does hold apart 
from rare exceptions. 

\medskip

There is a similar situation with respect to the location of the largest 
real zero, say $\beta$, of $L(s,\chi)$. Apart from rare exceptions one 
has the bound
\begin{equation}\label{eq:1.3}
1-\beta \gg  (\log D)^{-1}\ ,
\end{equation} 
and one expects that this always holds. Moreover, there is a close connection  
between the two phenomena and here we intend to study how close this is. In 
one direction this is clear, thanks to the following result of E. Hecke 
(see E. Landau [La]). 

\begin{theorem}\label{1} 
If (1.3) holds, then so does (1.2). 
\end{theorem}

\medskip

In this note we make a modest step toward the reversal of this 
implication. It is easy to see from the equation
$$
L'(s,\chi)= - \sum_{n\le D}\chi (n)(\log n)n^{-s} +O\bigl(D^{-\frac12}(\log D)^2\bigr)\ ,
$$
valid if $s-1\ll (\log D)^{-1}$, that
\begin{equation}\label{eq:1.4}
L'(s,\chi)\ll (\log D)^2 
\end{equation} 
there, and hence, by the mean-value theorem of differential calculus, 
that 
\begin{equation}\label{eq:1.5}
1-\beta \ll (\log D)^{-3} 
\end{equation} 
implies
\begin{equation}\label{eq:1.6}
L(1,\chi) \ll (\log D)^{-1}\ . 
\end{equation} 
Using deeper arguments we can say a bit more. 
\begin{theorem}\label{2} 
If $L(s,\chi)$ has a real zero $\beta$ with 
\begin{equation}\label{eq:1.7}
1-\beta \ll (\log D)^{-3}\, \log\log D\ , 
\end{equation} 
then (1.6) holds. 
\end{theorem}

\medskip

It seems surprising not to be able to do better. 
One expects, in the case of an extremely small value of $L(1,\chi)$,  
that $\chi$ mimics the M\"obius function and in such a situation 
$L'(s,\chi)$,  rather than being limited by (1.4), should be almost 
(though not quite) bounded, hence  
that (1.2) should imply something only slightly weaker than (1.3). 
The limitation in the bound of Theorem 2 comes from our 
imperfect knowledge about the complex zeros.
At the end of this note we describe (in Theorem 3) how the bound 
can be substantially sharpened if we assume the Riemann Hypothsis holds 
apart from an exceptional real zero. 

\medskip

There is an extensive literature on the subject discussed in this note. 
We encourage the reader to obtain a broader perspective through the 
publications [GoSc], [Go], [Pi], [MV], [GrSo], [SaZa]. 

\section {\bf Relations between $\beta$ and $L(1,\chi)$}
 
From now on we assume that $1-\beta \le (3\log D)^{-1}$. Denote 
$\lambda = 1 * \chi$ and 
$$
Z(s) =\zeta(s)L(s,\chi) = \sum_1^{\infty} \lambda (n)n^{-s}\ . 
$$
We evaluate the smoothly cropped sum 
\begin{equation}\label{eq:2.1}
S(x) =  \sum_{n\le x}\lambda (n)\bigl(1-\frac{n}{x}\bigr)n^{-\beta} 
\end{equation}  
by contour integration of $Z(s)$ as follows: 
\begin{equation*}
\begin{aligned}
S(x)  & = \frac{1}{2\pi i}\int_{(1)}Z(s+\beta)\frac{x^s}{s(s+1)}ds   \\
& = L(1,\chi)\, \frac{x^{1-\beta}}{(1-\beta)(2-\beta)} 
+ O\Bigl(x^{-1}D^{\frac12}\log D\Bigr)\ .
\end{aligned}
\end{equation*} 
Hence, 
\begin{prop}
Assume $L(s,\chi)$ has a real zero $\beta$ with 
 $1-\beta \le (3\log D)^{-1}$. Then 
\begin{equation}\label{eq:2.2}
L(1,\chi) \asymp (1-\beta)S(D)\ .
\end{equation}  
\end{prop} 

\medskip

Therefore, our problem reduces to the estimation of $S(D)$.  
Note that trivially $S(D) \gg 1$ and this gives Hecke's Theorem 1.

\section {\bf Upper Bound for $S(D)$}

We have 
\begin{equation}\label{eq:3.1}
S(D)\ll \sum_{n\le D}\frac{\lambda(n)}{n}
\ll \prod_{p\le D}\Bigl(1+\frac{\lambda(p)}{p}\Bigr)\ .
\end{equation} 
We split this product into 
\begin{equation}\label{eq:3.2}
\prod_{p\le B}\Bigl(1+\frac{\lambda(p)}{p}\Bigr)\, 
\le \prod_{p\le B}\Bigl(1+\frac{2}{p}\Bigr)\, \asymp (\log B)^2 
\end{equation} 
and 
\begin{equation}\label{eq:3.3}
\prod_{B<p\le D}\Bigl(1+\frac{\lambda(p)}{p}\Bigr) 
\le \exp \bigl(T(D)/\log B\bigr)\ ,
\end{equation} 
where 
\begin{equation}\label{eq:3.4}
T(D)=\sum_{p\le D}\frac{\lambda (p)}{p}\log p\ .
\end{equation} 

\medskip

Up to now our estimates have been rather simple, but to estimate 
$T(D)$ we require deeper tools, namely the Deuring-Heilbronn repulsion 
property of $\beta$. Put 
\begin{equation}\label{eq:3.5}
\eta = 1/(1-\beta)\log D \ge 3\ .
\end{equation} 
The repulsion property asserts that $L(s,\chi)$ has no zeros other than 
$\beta$ in the region
\begin{equation*}
\sigma > 1- \frac {c\log \eta}{\log D(|t|+1)}\ , \quad s = \sigma +it
\end{equation*} 
where $c$ is an absolute positive constant (cf. Th\'eor\`eme 16 
of E. Bombieri [B]). 

\medskip

A very nice way of using the repulsion property, together with a 
quite delicate estimate for the number of zeros of $L(s,\chi)$ in small 
discs centered on ${\rm Re}\, s = 1$ (cf. Ch. X, Lemma 2.1 of [Pr]), to bound 
sums of $\lambda (p)$ over primes has been given by D.R. Heath-Brown. We 
borrow from his work the estimate (Lemma 3 of [H-B1]) 
\begin{equation}\label{eq:3.6}
T(D)\ll (\log \eta)^{-\frac12}\log D\ . 
\end{equation} 
Choosing $\log B = (\log \eta)^{-\frac12}\log D$ we obtain 
\begin{equation}\label{eq:3.7}
S(D)\ll (\log \eta)^{-1}(\log D)^2\ . 
\end{equation} 

\section {\bf Conclusion}

From (2.2) and (3.7) we  get
\begin{equation}\label{eq:4.1}
L(1,\chi) \ll\frac{1-\beta} {\log \eta}(\log D)^2\ . 
\end{equation} 
In particular, if $\beta$ satisfies (1.7) then $\eta \gg \log D$ and (4.1) 
gives (1.6).

\section {\bf Remarks Behind the Scenes}

We take this opportunity to share some of our impressions about the 
nature of the arguments used in this note. The main issue is the question 
of how the exceptional zero is connected with the rarity of small primes 
which split in the quadratic field $\qb (\sqrt{\chi(-1)D})$. A quick 
connection is displayed in the bound (see (24.20) of [FI4]) 
\begin{equation}\label{eq:5.1}
\sum_{z<p\le x}\lambda (p) p^{-1} \ll (1-\beta) \log x\ , 
\end{equation} 
valid for $x>z\ge D^2$. This shows that if 
\begin{equation}\label{eq:5.2}
1-\beta =o(1/ \log D)\ , 
\end{equation} 
then the splitting primes in the segments $D^2<p\le D^A$ are very rare. 
One can easily deduce the same conclusion from the assumption 
(see (24.19) of [FI4])
\begin{equation}\label{eq:5.3}
L(1,\chi) = o (1/ \log D)\ . 
\end{equation} 

\medskip

The deficiency of such primes is the driving force for finding prime 
numbers in many interesting sequences; cf. [H-B1], [H-B2], [FI3]. One of 
these is 
the proof by Heath-Brown that the existence of infinitely many exceptional 
zeros implies the existence of infinitely many twin primes. Another example 
is the implication to primes of the form $p=a^2+b^6$. 

In our series of papers [FI1--FI5] we used assumptions of type 
(5.3) rather than (5.2) and for those applications they serve 
the same purpose. 

\medskip

Note that (5.1) says nothing about splitting primes which are very small 
relative to the conductor $D$. In all of these applications the rarity of 
small splitting primes was not needed, but for Theorem 2 it is essential. 
The current technology allows one to penetrate this territory, but only 
barely, due to the Deuring-Heilbronn repulsion property of the exceptional 
zero $\beta$. Heath-Brown's Lemma 3 of [H-B1] does just that! 

\medskip 

We include, for curiosity, an alternative derivation of (2.2). We 
consider the function
\begin{equation}\label{eq:5.4}
F(s) = (e^{\gamma}D)^sL(s+1, \chi)\ , 
\end{equation} 
which satisfies the conditions $F(\beta -1) = 0$, $F(0)=L(1,\chi)$,
\begin{equation*}
 F'(0)= L(1, \chi) (\log D +\gamma) 
+L'(1,\chi)\ 
\end{equation*} 
and $F''(s)\ll (\log D)^3$ if $s-1 \ll (\log D)^{-1}$. Hence, by the Taylor 
expansion of $F(\beta -1)$ at $s=0$ we get 
\begin{equation}\label{eq:5.5}
\begin{aligned}
L(1, \chi) & = (1-\beta)\bigl[L(1,\chi)(\log D + \gamma) +L'(1,\chi)\bigr]\\ 
& \quad\quad +O\bigl((1-\beta)^2(\log D)^3\bigr)\ .
\end{aligned}
\end{equation} 
On the other hand, we have (cf. (22.109) of [IK])
\begin{equation}\label{eq:5.6}
\sum_{n\le D} \lambda (n)n^{-1} = L(1,\chi)(\log D + \gamma) +L'(1,\chi) 
+O\bigl(D^{-\frac14}\log D\bigr)\ .
\end{equation} 
Combining (5.5) and (5.6) we obtain
\begin{equation}\label{eq:5.7}
L(1, \chi) = (1-\beta)\sum_{n\le D} \lambda (n)n^{-1}
+O\bigl((1-\beta)(\log D)^3(1-\beta +D^{-\frac14})\bigr)\ .
\end{equation} 
For our purposes, (5.7) and (2.2) amount to the same thing. 

\medskip

Note that, under the assumption (5.3) the formula (5.7) implies 
\begin{equation}\label{eq:5.8}
L'(1, \chi) \sim \sum_{n\le D} \lambda (n)n^{-1}\ ,
\end{equation} 
which can be compared with (22.117) of [IK]. 

\medskip

There are infinitely many real primitive characters $\chi (\mod D)$ 
with $D$ prime and $\chi (p) =1$ for every $p<c\log D$. For such characters 
we have
\begin{equation}\label{eq:5.9}
\sum_{n\le D} \lambda (n)n^{-1}\gg (\log\log D)^2\ .
\end{equation} 
Hence, if the largest real zero of $L(s,\chi)$ satisfies (3.5), we have 
\begin{equation}\label{eq:5.10}
L(1,\chi) \gg (1-\beta)(\log\log D)^2\ .
\end{equation} 
Therefore, (1.6) cannot hold for such special discriminants unless  
$L(s,\chi)$ has a zero $\beta$ with 
$1-\beta \ll (\log D)^{-1} (\log\log D)^{-2}$. 

\section {\bf The Ultimate Deuring-Heilbronn Phenomenon}

In this final section we investigate the extent of improvement in the 
conclusion of Theorem 2 
which could be obtained if one assumed the GRH for $L(s,\chi)$ apart 
from one real zero $\beta$. Let us assume that $\beta >\frac34$ 
is the only zero of $L(s,\chi)$ in ${\rm Re}\, s >\frac34$. We proceed along 
the lines of J. E. Littlewood [Li], beginning with the formula 
(cf. (5.58) of [IK])
\begin{equation}\label{eq:6.1}
\begin{aligned}
-\frac{L'}{L}(\sigma,\chi)& =\sum_{n\le x}\chi(n)\Lambda(n)
\bigl(1-\frac{n}{x}\bigr)n^{-\sigma}\\ & -\sum_{\rho}(\rho -\sigma)^{-1}
(\rho -\sigma +1)^{-1}x^{\rho -\sigma} +O\bigl(x^{-\frac12} \log D\bigr)
\end{aligned}
\end{equation} 
which is valid for $1\le \sigma\le 5/4$ with any $x\ge 1$, the implied 
constant being absolute. Here $\rho$ runs through the zeros of $L(s,\chi)$ 
in the critical strip. Separating $\rho = \beta$ and estimating the other 
terms trivially we find
\begin{equation}\label{eq:6.2}
\begin{aligned}
-\frac{L'}{L}(\sigma,\chi)& =\sum_{n\le x}\chi(n)\Lambda(n)n^{-\sigma} \\ & 
-(\beta - \sigma)^{-1}
(\beta - \sigma +1)^{-1}x^{\beta -\sigma} +O\bigl(1+ x^{-\frac14} \log D\bigr)\ . 
\end{aligned}
\end{equation} 
We take $x=(\log D)^4$ so the error term in (6.2) is bounded. Integrating (6.2) 
over $1\le \sigma \le 5/4$ we obtain 
\begin{equation*}
\log L(1,\chi) =\sum_{p\le x}\frac{\chi(p)}{p}- 
\int_{1-\beta}^{5/4 -\beta}x^{-t}t^{-1}(1-t)^{-1}dt +O(1)\ . 
\end{equation*} 
Up to a bounded error term, the integral is equal to
\begin{equation*}
\int_{1-\beta}^{\infty}x^{-t}t^{-1}dt =  \int_{\vt}^{\infty}e^{-t}t^{-1}dt 
= - e^{-\vt}\log {\vt} + \int_{\vt}^{\infty} e^{-t}(\log t) dt\ , 
\end{equation*} 
where $\vt =(1-\beta) \log x$. Here, we have 
$e^{-\vt}\log {\vt} = (\log \vt)/ (\vt+1) +O(1)$ and the last integral is 
bounded. Hence, 
\begin{equation*}
\log L(1,\chi) =\sum_{p\le x}\chi(p)p^{-1}  
+(\log{\vt})/(\vt + 1) +O(1)\ . 
\end{equation*} 
In other words 
\begin{equation}\label{eq:6.3}
L(1,\chi)\asymp \omega \exp\Bigl(\sum_{p\le \log D}\chi(p)p^{-1} \bigr)
\end{equation} 
where 
\begin{equation}\label{eq:6.4}
\omega =\min\{1,(1-\beta)\log\log D\}\ . 
\end{equation} 
This proves 
\begin{theorem}\label{3} 
Assuming that $\beta >3/4$ is the only zero of $L(s,\chi)$ in 
${\rm Re}\, s > 3/4$ 
we have (6.3). Hence
\begin{equation}\label{eq:6.5}
\omega (\log\log D)^{-1} \ll L(1,\chi) \ll \omega (\log\log D)\ . 
\end{equation} 
In particular, if $1-\beta \ll (\log\log D)^{-1}$, then 
\begin{equation}\label{eq:6.6}
1-\beta \ll L(1,\chi) \ll (1-\beta) (\log\log D)^2\ . 
\end{equation} 
\end{theorem}

\medskip 
Department of Mathematics, University of Toronto

Toronto, Ontario M5S 2E4, Canada 

\medskip

Department of Mathematics, Rutgers University

Piscataway, NJ 08903, USA

\end{document}